\documentclass{amsart}
\usepackage{amssymb,latexsym,enumerate,url}

\newtheorem{Thm}[equation]{Theorem}
\newtheorem{Prop}[equation]{Proposition}
\newtheorem{Lem}[equation]{Lemma}
\newtheorem{Cor}[equation]{Corollary}
\newtheorem*{Cor*}{Corollary}

\theoremstyle{remark}

\newtheorem*{Rem*}{Remark}
\theoremstyle{definition}

\newtheorem*{Not*}{Notation}

\numberwithin{equation}{section}
\DeclareMathOperator{\im}{im}
\DeclareMathOperator{\FP}{FP}
\DeclareMathOperator{\Hom}{Hom}
\DeclareMathOperator{\Diff}{Diff}
\begin{document}

\date{%
Tue May 11 07:51:40 EDT 2010}

\title{Dimensions of $\ell^p$-cohomology groups}

\author[M. S. Grinshpon]{Mark S. Grinshpon}
\address{Department of Mathematics and Statistics\\
Georgia State University\\
Atlanta, GA 30303\\
USA}
\email{matmsg@langate.gsu.edu}
\urladdr{\url{http://www2.gsu.edu/~wwwmat/faculty_staff/instructors/grinshpon.html}}

\author[P. A. Linnell]{Peter A. Linnell}
\address{Department of Mathematics\\
Virginia Tech\\
Blacksburg\\
VA 24061-0123\\
USA}
\email{plinnell@math.vt.edu}
\urladdr{\url{http://www.math.vt.edu/people/plinnell/}}

\author[M. J. Puls]{Michael J. Puls}
\address{Department of Mathematics, John Jay College--CUNY\\
445 West 59th Street, New York\\
NY 10019, USA}
\email{mpuls@jjay.cuny.edu}
\thanks{The research of the third author was partially
supported by PSC-CUNY grant 62598-00 40}

\begin{abstract}
Let $G$ be an infinite discrete group of type
$\FP_{\infty}$ and let $1 < p \in
\mathbb{R}$.  We prove that the $\ell^p$-homology
and cohomology groups of $G$ are either 0 or infinite
dimensional.  We also show that the cardinality of the $p$-harmonic
boundary of a finitely generated group is either 0, 1, or $\infty$.
\end{abstract}

\keywords{projective resolution, group of type $\FP_n$,
$\ell^p$-cohomology, translation invariant functional}

\subjclass[2000]{Primary: 43A15; Secondary: 20J06}

\maketitle

\section{Introduction}

Let $G$ be a discrete group and for $1\le p \in \mathbb{R}$, let
\[
\ell^p(G) =\{ \sum_{x\in G} a_xx \mid a_x\in \mathbb{C} \text{ and }
\sum_{x \in G} |a_x|^p < \infty\}.
\]
This is a complex Banach space with
respect to the norm $\| f\|_p := \bigl(\sum_{x\in G} |a_x|^p
\bigr)^{1/p}$.  Also $G$ acts on the left of $\ell^p(G)$ according to
the rule $g \sum_{x\in G} a_xx = \sum_{x\in G} a_xgx$, and similarly
on the right according to the rule $\bigl( \sum_{x\in G} a_xx \bigr)g
= \sum_{x\in G} a_xxg$.  These actions make $\ell^p(G)$ into a
$\mathbb{C}G$-bimodule.
Suppose we are given a free resolution of
the trivial $\mathbb{C}G$-module $\mathbb{C}$ with free
right $\mathbb{C}G$-modules:
\begin{equation} \label{Esequence}
\cdots \longrightarrow
\mathbb{C}G^{e_{n+1}} \overset{d_n}{\longrightarrow}
\mathbb{C}G^{e_n}
\overset{d_{n-1}}{\longrightarrow}
\cdots \overset{d_1}{\longrightarrow}
\mathbb{C}G^{e_1} \overset{d_0}{\longrightarrow}
\mathbb{C}G \longrightarrow
\mathbb{C} \longrightarrow 0.
\end{equation}
Let
\begin{align}
d_n^* &\colon \Hom_{\mathbb{C}G} (\mathbb{C}G^{e_n},\ell^p(G))
\longrightarrow
\Hom_{\mathbb{C}G}(\mathbb{C}G^{e_{n+1}},\ell^p(G)),\label{Ecohom}\\
d^n_* &\colon \mathbb{C}G^{e_{n+1}} \otimes_{\mathbb{C}G} \ell^p(G)
\longrightarrow
\mathbb{C}G^{e_n} \otimes_{\mathbb{C}G} \ell^p(G) \label{Ehom}
\end{align}
be the maps induced by $d_n$; for convenience we write
$d_{-1}^* = d^{-1}_* = 0$.
Then one has the usual (unreduced) cohomology and homology groups
\begin{align*}
H^n(G,\ell^p(G)) &= \ker d_n^*/\im d_{n-1}^*, \\
H_n(G,\ell^p(G)) &= \ker d_*^{n-1}/\im d_*^n.
\end{align*}
We will be interested in the case $e_n < \infty$ for all $n$, whence
$\Hom_{\mathbb{C}G} (\mathbb{C}G^{e_n},\ell^p(G)) \cong
\ell^p(G)^{e_n}$ as left $\mathbb{C}G$-modules and
$\mathbb{C}G^{e_n} \otimes_{\mathbb{C}G} \ell^p(G) \cong
\ell^p(G)^{e_n}$ as right $\mathbb{C}G$-modules.
(Recall that $\Hom_{\mathbb{C}G}(\mathbb{C}G^{e_n}, \ell^p(G))$
consists of right $\mathbb{C}G$-maps $\theta \colon
\mathbb{C}G^{e_n} \to \ell^p(G)$ with left $G$-action defined by
$(g\theta)(\alpha) = g(\theta \alpha)$ for $\alpha \in
\mathbb{C}G^{e_n}$.  Also $\mathbb{C}G^{e_{n+1}}
\otimes_{\mathbb{C}G} \ell^p(G)$ is the tensor product of the right
$\mathbb{C}G$-module $\mathbb{C}G^{e_{n+1}}$ with the left
$\mathbb{C}G$-module $\ell^p(G)$ with right $G$-action defined by
$(\alpha \otimes u)g = \alpha \otimes ug$.)

Bekka and Valette \cite[Corollary 8]{BekkaValette97} proved that if
$G$ is a finitely generated group, then $H^1(G,\ell^2(G))$
is either zero or infinite dimensional.  The motivation behind this
paper is to see if this result holds for arbitrary $1 < p \in
\mathbb{R}$.  It turns out that this result remains true not only for
$1 < p \in \mathbb{R}$, but also for the other homology and
cohomology groups.  Recall that $G$ is of type $\FP_n$ over
$\mathbb{C}$ if there exists a resolution \eqref{Esequence}
which has $e_d$ finite for all $d \le n$, and $G$ is of type
$\FP_{\infty}$ if it is of type $\FP_n$ for all $n\in \mathbb{N}$.
Furthermore $G$ is of type $\FP_1$ over $\mathbb{C}$ if
and only if $G$ is finitely generated.  We shall prove
\begin{Thm} \label{Tmain1}
Let $d,n$ be non-negative integers and let $G$ be an infinite
group of type
$\FP_n$ over $\mathbb{C}$.  Let $1 < p \in \mathbb{R}$.  Then
\begin{enumerate}[\normalfont(i)]
\item $H^d(G,\ell^p(G))$ is either 0 or has infinite
$\mathbb{C}$-dimension for all $d \le n$. \label{T1}
\item $H_d(G,\ell^p(G))$ is either 0 or has infinite
$\mathbb{C}$-dimension for all $d \le n$. \label{T2}
\end{enumerate}
\end{Thm}
Theorem \ref{Tmain1} immediately yields
\begin{Cor*}
Let $G$ be an infinite group of type $\FP_{\infty}$
over $\mathbb{C}$ and let $d$ be a non-negative integer.
Let $1 < p \in \mathbb{R}$.  Then
\begin{enumerate}[\normalfont(i)]
\item $H^d(G,\ell^p(G))$ is either 0 or has infinite
$\mathbb{C}$-dimension.
\item $H_d(G,\ell^p(G))$ is either 0 or has infinite
$\mathbb{C}$-dimension.
\end{enumerate}
\end{Cor*}
We deduce Theorem \ref{Tmain1} from our main theorem:
\begin{Thm} \label{Tmain2}
Let $G$ be an infinite group, let $m$ be a non-negative integer,
and let $A \subseteq B$ be closed left $G$-invariant
subspaces of $\ell^p(G)^m$.
Then either $A=B$ or $B/A$ has infinite dimension over $\mathbb{C}$.
\end{Thm}
Of course, Theorem \ref{Tmain2} remains true if we replace ``left"
with ``right".

The layout of this paper is as follows.  In Section
\ref{SPreliminaries} we give some definitions and recall some
well-known results.  In Section \ref{Smain} we prove Theorems
\ref{Tmain1} and \ref{Tmain2}.
In Section \ref{Stilf}, we shall use Theorem \ref{Tmain2} to obtain
a result concerning the cardinality of the $p$-harmonic boundary of a
finitely generated group and to prove a result about translation
invariant functionals on a certain function space of functions on a
finitely generated group.

\section{Preliminaries} \label{SPreliminaries}

We denote the positive integers by $\mathbb{N}$.  Let $1 \le p
\in \mathbb{R}$.  Then for $\alpha = \sum_{g\in G}
\alpha_gg \in \ell^1(G)$ and
$\beta = \sum_{g\in G} \beta_gg \in \ell^p(G)$, we define
convolution by
\[
\alpha \beta = \sum_{g,h \in G} \alpha_g\beta_h gh = \sum_{g\in G}
\bigl( \sum_{x \in G} \alpha_{gx^{-1}}\beta_x \bigr) g \in \ell^p(G).
\]
Young's inequality \cite[32D]{Loomis53} tells us that
\begin{equation} \label{EYoung}
\|\alpha\beta\|_p \le \|\alpha\|_1 \|\beta\|_p.
\end{equation}
Thus in particular $\ell^1(G)$ is a ring with multiplication being
convolution.  Let $m$ be a
non-negative integer.  While different norms can be defined on a
finite direct sum of normed spaces, they are all equivalent \cite[\S
1.8]{Megginson98}.  The most natural and consistent choice for
$\ell^p(G)^m$ is using the $p$-norm:
\[
\|(u_1,\dots,u_m)\|_p = \bigl( \sum_{k=1}^m \| u_k \|_p^p\bigr)^{1/p}
\text{ for all } (u_1,\dots,u_m) \in \ell^p(G)^m.
\]
The inequality \eqref{EYoung} still holds for $u \in \ell^1(G)$
and $v = (v_1,\dots,v_m)
\in \ell^p(G)^m$ with convolution defined componentwise, because
\begin{equation} \label{EYoung1}
\begin{aligned}
\|uv\|_p^p = \|u(v_1,\dots,v_m)\|_p^p &=\sum_{k=1}^m \|uv_k\|_p^p
\le \sum_{k=1}^m \bigl(\|u\|_1 \|v_k\|_p)^p\\
&= \sum_{k=1}^m \|u\|_1^p \|v_k\|_p^p = \|u\|_1^p \sum_{k=1}^m
\|v_k\|_p^p = \|u\|_1^p \|v\|_p^p.
\end{aligned}
\end{equation}
Similarly for $u \in \ell^p(G)$ and $v = (v_1,\dots, v_m)
\in \ell^1(G)^m$, we have
\begin{equation} \label{EYoung2}
\| uv \|_p^p \le \|u\|_p^p \|v\|_1^p.
\end{equation}
Note that \eqref{EYoung1} tells us that $\ell^p(G)^m$ is a left
$\ell^1(G)$-module, and that closed left $G$-invariant subspaces of
$\ell^p(G)^m$ are left $\ell^1(G)$-submodules.

We shall also need the following two well-known results.
\begin{Lem} \label{Ldense}
Suppose $n$ is a positive integer and for $1 \le k \le n$
we have bounded linear operators $T_k \colon B\to B$ on a
normed space $B$ such that the range $T_k(B)$ is dense in $B$ for
each $T_k$.  Then the range of $T_1 \cdots T_n$ is also dense in $B$.
\end{Lem}
\begin{proof}
First we prove the claim for $n=2$.  Let $\epsilon > 0$ and $b \in
B$.  There exists $b_1 \in B$ such that $\| b-T_1 b_1 \| <
\epsilon/2$, and then there exists $b_2 \in B$ such that $\|b_1 -
T_2b_2\| \le \frac{\epsilon}{2\|T_1\|}$.  Thus
\begin{align*}
\| b - T_1T_2b_2 \| &\le \|b-T_1b_1\| + \| T_1b_1 - T_1T_2b_2 \| \\
&\le \| b-T_1b_1 \| + \| T_1\| \| b_1 -T_2b_2\| \\
&< \frac{\epsilon}{2} + \|T_1\| \cdot \frac{\epsilon}{2\| T_1 \|} =
\epsilon.
\end{align*}
The lemma now follows by induction on $n$.
\end{proof}

\begin{Lem} \label{Lclosed}
Let $T\colon A\to B$ be a bounded linear operator between the Banach
spaces $A$ and $B$.  If $T(A)$ has finite codimension in $B$, then
$T(A)$ is closed in $B$.
\end{Lem}
\begin{proof}
See \cite[p.~95, Exercise (1), \S3.4]{Arveson02}
\end{proof}

\section{Proof of the main theorems} \label{Smain}

The critical case in the proof of Theorem \ref{Tmain2} is when $G$ is
infinite cyclic, and the reader will understand most of the proof by
studying this special situation.  To prove the result in general,
we have had to repeat some arguments almost verbatim several
times.  However we have chosen to give full details over brevity and
clarity.

\begin{proof}[Proof of Theorem \ref{Tmain2}.]
We will assume that $B/A$ is finite dimensional and will prove that
$A = B$.  First suppose $G$ has an element $g$ of infinite order.
Write $H = \langle g\rangle$.  Note that for $\alpha \in \ell^1(G)$
and $\beta \in \ell^p(G)^m$, we have $\| \alpha \beta\|_p \le
\| \alpha\|_1 \| \| \beta \|_p$ by \eqref{EYoung1},
thus in particular $\ell^p(G)^m$ is a left $\ell^1(H)$-module,
and $A$ and $B$ are left $\ell^1(H)$-submodules.
The action of $g$ on the finite dimensional vector space $B/A$ has a
minimal polynomial, i.e.\ there exists $F(x) \in \mathbb{C}[x]$ such
that $F(g) = 0$ on $B/A$, and therefore $F(g)b \in A$ for all $b \in
B$.  Factor $F(x)$ into linear factors
and notice that if $|\omega| \ne 1$, then
$(g-\omega)$ is invertible in $\ell^1(H)$.  Thus since $A$ and $B$
are $\ell^1(H)$-invariant, we may assume that $F(g)$ consists of
factors $(g-\omega)$ with $|\omega| = 1$ only.  If we prove that
$F(g)B$ is dense in $B$, that will imply that $A=B$.

Fix $\omega$ with $|\omega|=1$ and
for $n \in \mathbb{N}$, let $x_n = \frac{1}{n} \sum_{k=1}^n
\omega^{-k}g^k$.  Note that
\begin{equation} \label{Enorm}
\| x_n \|_p = \left ( n \cdot \frac{1}{n^p} \right)^{1/p} =
n^{(1-p)/p},
\end{equation}
consequently $\lim_{n\to \infty} \|x_n\|_p = 0$.

Now pick arbitrary $b \in B$ and $\epsilon > 0$.  Since
$\mathbb{C}G^m$ is dense in $\ell^p(G)^m$, there exists
$c \in \mathbb{C}G^m$ such that $\| b-c \|_p < \epsilon/2$.  Then
we may choose $n \in \mathbb{N}$ such that $\| x_n \|_p <
\frac{\epsilon}{2 \| c\|_1}$ and we have
\begin{align*}
\|x_nb\|_p &= \|x_n(b-c) + x_nc\|_p \le \|x_n(b-c) \|_p +
\| x_nc\|_p \\
&\le \| x_n \|_1 \| b-c \|_p + \| c\|_1 \| x_n\|_p
\quad\text{by \eqref{EYoung1} and \eqref{EYoung2}} \\
&< 1 \cdot \frac{\epsilon}{2} + \|c\|_1 \cdot \frac{\epsilon}{2
\|c\|_1} = \epsilon.
\end{align*}
Thus $\| b - (1-x_n)b\|_p = \|x_n b\|_p < \epsilon$.  Now note that
the homomorphism $\mathbb{C}H \to \mathbb{C}$ induced
by the identity on $\mathbb{C}$ and sending $g$ to $\omega$ has
$(1-x_n)$ in its kernel, consequently we may write $1-x_n =
(g-\omega)d$, where $d \in \mathbb{C}G$,
and we deduce that $(1-x_n)b \in
(g-\omega)B$.  Thus $(g-\omega)B$ is dense in $B$.  Since the product
of operators with dense ranges has dense range by Lemma \ref{Ldense},
we conclude that $F(g)B$ is dense in $B$.

Therefore we may assume that every element of $G$ has finite order.
Let $N$ denote the kernel of the action of $G$ on $B/A$.  Suppose $N$
is infinite.  Choose an infinite sequence $\{g_1,g_2,\dots \}$ of
distinct elements of $N$ and let $x_n = \sum_{k=1}^n \frac{1}{n} g_k
\in \mathbb{C}G$.  Let $b \in B$, let $\epsilon > 0$ and follow the
argument above.  Since $\|x_n\|_p = n^{(1-p)/p}$, we see that
$\lim_{n\to \infty} \|x_n\|_p = 0$.  Also $\mathbb{C}G^m$
is dense $\ell^p(G)^m$, hence there exists $c \in
\mathbb{C}G^m$ such that $\|b-c\|_p < \epsilon/2$.  Then
we may choose $n \in \mathbb{N}$ such that $\| x_n \|_p <
\frac{\epsilon}{2 \| c\|_1}$ and we have
\begin{align*}
\|x_nb\|_p &= \|x_n(b-c) + x_nc\|_p \le \|x_n(b-c) \|_p +
\| x_nc\|_p \\
&\le \| x_n \|_1 \| b-c \|_p + \| c\|_1 \| x_n\|_p
\quad\text{by \eqref{EYoung1} and \eqref{EYoung2}} \\
&< 1 \cdot \frac{\epsilon}{2} + \|c\|_1 \cdot \frac{\epsilon}{2
\|c\|_1} = \epsilon.
\end{align*}
Thus $\| b - (1-x_n)b\|_p = \|x_n b\|_p < \epsilon$.  Since $(1-x_n)b
\in A$ for all $n \in \mathbb{N}$,
we see that $A$ is dense in $B$ and we conclude that $A=B$.

Therefore we may assume that $N$ is finite, so $G/N$ is an infinite
torsion group, and its action on $B/A$ tells us that it is also a
linear group over $\mathbb{C}$.  By a theorem of Schur
\cite[cf.~1.L.4]{KegelWehrfritz73}, there is a normal abelian
subgroup $K/N$ of finite index in $G/N$.  Since a simple
$\mathbb{C}[K/N]$-module has dimension one over $\mathbb{C}$, there
is a one-dimensional $K$-invariant subspace $U/A$ of $B/A$.
Again follow the proof above.  Choose an infinite sequence
$\{g_1,g_2,\dots \}$ of distinct elements of $K/N$.  Then there
exist $\omega_1,\omega_2,\ldots \in \mathbb{C}$ with $|\omega_i| =1$
such that $(g_i - \omega_i) U \subseteq A$.  As before for $n \in
\mathbb{N}$, set $x_n = \sum_{k=1}^n \frac{1}{n} \omega_k^{-1}g_k$.
Again (cf.~\eqref{Enorm}), $\|x_n\|_p = n^{(1-p)/p} \to 0$
as $n \to \infty$.

Now pick arbitrary $u \in U$ and $\epsilon > 0$.  Since
$\mathbb{C}G^m$ is dense in $\ell^p(G)^m$, there exists
$c \in \mathbb{C}G^m$ such that $\| u-c \|_p < \epsilon/2$.  Then
we may choose $n \in \mathbb{N}$ such that $\| x_n \|_p <
\frac{\epsilon}{2 \| c\|_1}$ and we have
\begin{align*}
\|x_nu\|_p &= \|x_n(u-c) + x_nc\|_p \le \|x_n(u-c) \|_p +
\| x_nc\|_p \\
&\le \| x_n \|_1 \| u-c \|_p + \| c\|_1 \| x_n\|_p
\quad\text{by \eqref{EYoung1} and \eqref{EYoung2}} \\
&< 1 \cdot \frac{\epsilon}{2} + \|c\|_1 \cdot \frac{\epsilon}{2
\|c\|_1} = \epsilon.
\end{align*}
Since $(1-x_n)u \in A$ for all $n$, we deduce that $A$ is dense in
$U$, a contradiction.  This completes the proof of Theorem
\ref{Tmain2}.
\end{proof}

\begin{proof}[Deduction of Theorem \ref{Tmain1} from Theorem
\ref{Tmain2}]
For Theorem \ref{Tmain1}\eqref{T1}, note that the maps $d_n^* \colon
\ell^p(G)^{e_n} \to \ell^p(G)^{e_{n+1}}$ are continuous, because they
are given by right multiplication by an $e_n \times e_{n+1}$ matrix
with entries in $\mathbb{C}G$.  Thus if $\im d_{n-1}^*$ has finite
codimension in $\ker d_n^*$, it will be closed by Lemma
\ref{Lclosed}.  Theorem \ref{Tmain1}\eqref{T1} now follows from
Theorem \ref{Tmain2}.  The proof of Theorem \ref{Tmain1}\eqref{T2} is
almost exactly the same, except we need to deal with right
$G$-invariant subspaces of $\ell^p(G)^{e_n}$.
\end{proof}

We can also prove results for the corresponding real Banach spaces.
For a group $G$ and $1\le p \in \mathbb{R}$, let
\[
l^p(G) =\{ \sum_{x\in G} a_xx \mid a_x\in \mathbb{R} \text{ and }
\sum_{x \in G} |a_x|^p < \infty\}.
\]
This is a real Banach space with
respect to the norm $\| f\|_p := \bigl(\sum_{x\in G} |a_x|^p
\bigr)^{1/p}$.  Also $G$ acts on the left of $l^p(G)$ according to
the rule $g \sum_{x\in G} a_xx = \sum_{x\in G} a_xgx$.  Then we have
\begin{Cor} \label{Cmain2}
Let $G$ be an infinite group, let $m$ be a non-negative integer,
and let $A \subseteq B$ be closed left $G$-invariant
subspaces of $l^p(G)^m$.
Then either $A=B$ or $B/A$ has infinite dimension over $\mathbb{R}$.
\end{Cor}
\begin{proof}
We can regard $A$ and $B$ as closed $G$-invariant real subspaces of
$\ell^p(G)^m$.  Set $X = A + iA$ and $Y = B + iB$.  Then $X$ and
$Y$ are closed left $G$-invariant complex subspaces of $\ell^p(G)^m$.
Since either $Y= X$ or $Y/X$ has infinite $\mathbb{C}$-dimension by
Theorem \ref{Tmain2}, we see that either $B=A$ or $B/A$ has infinite
$\mathbb{R}$-dimension and the result follows.
\end{proof}

\section{Applications to finitely generated groups} \label{Stilf}

In this section we will use Corollary \ref{Cmain2} to obtain some new
results concerning finitely generated infinite groups.  Let
$\mathcal{F}(G)$ denote the set of all real valued functions on $G$.
This has a left and right $G$-action given by $(gf)(x) = f(g^{-1}x)$
and $(fg)(x) = f(xg^{-1})$ for $f \in \mathcal{F}(G)$ and $g,x \in
G$, respectively.  We will view $l^p(G)$ as
$\{f \in \mathcal{F}(G) \mid \sum_{x\in G} |f(x)|^p < \infty\}$.
To make this identification, we
send $f$ to $\sum_{x\in G} f(x)x$.  Also $l^{\infty}(G) = \{f \in
\mathcal{F}(G) \mid \sup_{x\in G} |f(x)| < \infty\}$ with norm $\| f
\|_{\infty} = \sup_{x\in G} |f(x)|$.  Finally $\mathbb{R}G$ will
denote the functions in $\mathcal{F}(G)$ with finite support.

Throughout this section, $p$ will always denote a real number
greater than one and $G$ will be a group with a finite symmetric
generating set $S$ (so $S = S^{-1}$).  For a
real-valued function $f$ on $G$, we define the $p$-th
power of the \emph{gradient}, the \emph{$p$-Dirichlet sum}, and the
\emph{$p$-Laplacian} of $g \in G$ by
\begin{align*}
| Df(g)|^p &= \sum_{s\in S} |f(g) - f(gs)|^p,\\
I_p(f) &= \sum_{g \in G} |Df(g)|^p, \text{ and }\\
\Delta_p f(g) &= \sum_{s\in S} |f(gs) - f(g)|^{p-2} (f(gs)
- f(g)),
\end{align*}
respectively.
In the case $1 < p < 2$, we make the convention that $|f(gs) -
f(g)|^{p-2} (f(gs) - f(g)) = 0$ if $f(gs) = f(g)$.  We shall
say that $f$ is \emph{$p$-Dirichlet finite} if $I_p(f) < \infty$.
The set of all $p$-Dirichlet finite functions on $G$ will be denoted
by $D_p(G)$.  A function $f$ is said to be $p$-harmonic if $\Delta_p
f(g) = 0$ for all $g \in G$.  The set $HD_p(G)$ will consist of the
$p$-harmonic functions contained in $D_p(G)$.  We identify the
constant functions on $G$ with $\mathbb{R}$.  Observe that
$\mathbb{R}$ is contained in $HD_p(G)$.  Endowed with the norm
\[
\| f\|_{D_p} := (I_p(f) + | f(e)|^p)^{1/p},
\]
$D_p(G)$ is a reflexive Banach space, where $e$ is the identity
element of $G$ and $f \in D_p(G)$.  For $X \subseteq D_p(G)$, let
$\overline{X}_{D_p}$ indicate its closure in the $D_p$-norm and let
$B(X)$ denote the bounded functions in $X$, that is $X\cap
l^{\infty}(G)$; sometimes we will write $BX$ for $B(X)$.  The set
$BD_p(G)$ is closed under the usual operations of scalar
multiplication and addition.  Also $BD_p(G)$ is a reflexive Banach
space under the norm
\[
\| f\|_{BD_p} := (I_p(f))^{1/p} + \|f\|_{\infty},
\]
where $f \in BD_p(G)$.  Furthermore, $\| fh\|_{BD_p} \le \|f\|_{BD_p}
\|h \|_{BD_p}$ for $f,h \in BD_p(G)$.  Thus $BD_p(G)$ is an abelian
Banach algebra.  For $Y \subseteq BD_p(G)$, let $\overline{Y}_{BD_p}$
denote its closure in the $BD_p$-norm.  Note that if $f \in BD_p(G)$,
then $\| f\|_{D_p} \le \|f\|_{BD_p}$ and that
$B(\overline{\mathbb{R}G}_{D_p}) = B(\overline{l^p(G)}_{D_p})$
is a closed ideal in $BD_p(G)$.

Our first application of Corollary \ref{Cmain2} will be concerned with
the cardinality of the $p$-harmonic boundary of $G$, which we now
define.  For a more detailed discussion of this boundary, see
\cite{Puls08}.  Let $Sp(BD_p(G))$ denote the set of complex-valued
characters on $BD_p(G)$, that is nonzero $*$-homomorphisms from
$BD_p(G)$ to $\mathbb{C}$.  Then with respect to the weak
$*$-topology, $Sp(BD_p(G))$ is a compact Hausdorff space.  Given a
topological space $X$, let $C(X)$ denote the ring of continuous
complex-valued functions on $X$.  The Gelfand transform defined by
$\hat{f}(\chi) = \chi(f)$ yields a monomorphism of Banach algebras
from $BD_p(G)$ into $C(Sp(BD_p(G)))$ with dense image.  Furthermore
the map $\iota \colon G \to Sp(BD_p(G))$ given by $(\iota(g))(f) =
f(g)$ is an injection, and $\iota(G)$ is an open dense subset of
$Sp(BD_p(G))$.  The \emph{$p$-Royden boundary} of $G$, which we shall
denote by $R_p(G)$, is the compact set $Sp(BD_p(G)) \setminus
\iota(G)$.  The \emph{$p$-harmonic boundary} of $G$ is the following
subset of $R_p(G)$:
\[
\partial_p(G) := \{\chi \in R_p(G) \mid \hat{f}(\chi) = 0 \text{ for
all } f \in B(\overline{\mathbb{R}G}_{D_p})\}.
\]
We can now state
\begin{Thm}
Let $1 < p \in \mathbb{R}$ and let $G$ be a finitely generated
infinite group.  Then the cardinality of $\partial_p(G)$ is
either 0, 1 or $\infty$.
\end{Thm}
\begin{proof}
Let $S := \{s_1,\dots,s_d\}$ be a symmetric generating set for $G$.
We will use the results of \cite{Puls08} with $G$ being the Cayley
graph of $G$ with respect to the generating set $S$; thus the
vertices of this graph are the elements of $G$, and $g_1,g_2 \in G$
are joined by an edge if and only if $g_1 = g_2s$ for some $s\in S$.
If $1 \in B(\overline{\mathbb{R}G}_{D_p})$, then $\partial_p(G) =
\emptyset$ by \cite[Theorem 2.1 and Proposition 4.2]{Puls08}.
Thus we will assume that $1
\notin B(\overline{\mathbb{R}G}_{D_p})$.  If $BHD_p(G) =
\mathbb{R}$, then \cite[Theorem 4.11]{Puls08} says that
$|\partial_p(G)| = 1$.  Now suppose $|\partial_p(G)| > 1$.  We will
complete the proof of the theorem by showing that if
$|\partial_p(G)|$ is finite, then there exist two closed left
$G$-invariant subspaces $A$ and $B$ of $l^p(G)^d$
that violate Corollary \ref{Cmain2}.  We start by
showing how $D_p(G)$ is related to $l^p(G)^d$.  Define a continuous
linear map $D_p(G) \to l^p(G)^d$ by $\theta(f) =
(f(s_1-1),\dots,f(s_d-1))$.  Then $\ker\theta = \mathbb{R}$ and so
$\theta$ induces an embedding $\theta' \colon D_p(G)/\mathbb{R}
\hookrightarrow l^p(G)^d$.  Clearly $D_p(G)/\mathbb{R}$ is a Banach
space under the norm induced by the norm on $D_p(G)$, and $\theta'$
preserves this norm and also the left $G$-action.

We now construct the subspaces of $l^p(G)^d$ that will give us our
contradiction.  The Gelfand transform yields a homomorphism of
$BD_p(G)/B(\overline{\mathbb{R}G}_{D_p})$ onto a dense subspace of
$C(\partial_p(G))$.  Now \cite[Theorem 4.9]{Puls08} shows that if $f
\in BD_p(G)$, then $f \in B(\overline{\mathbb{R}G}_{D_p})$ if and only
if $\hat{f} = 0$ on $\partial_p(G)$.  This tells us that if
$|\partial_p(G)| < \infty$, then
$\dim_{\mathbb{R}}\bigr(BD_p(G)/B(\overline{\mathbb{R}G}_{D_p})\bigl)
= |\partial_p(G)|$.  Suppose $1 < | \partial_p(G) | < \infty$.  Then
\begin{align*}
&1\le \dim_{\mathbb{R}} \bigl(BD_p(G)/(B(\overline{\mathbb{R}G}_{D_p})
\oplus \mathbb{R})\bigr) < \infty\\
\intertext{and we deduce that}
&1\le \dim_{\mathbb{R}} \bigl((BD_p(G) +
\overline{\mathbb{R}G}_{D_p})/\overline{\mathbb{R}\oplus \mathbb
{R}G}_{D_p}\bigr) < \infty.
\end{align*}
It follows that $\theta'\bigl((BD_p(G)) +
\overline{\mathbb{R}G}_{D_p}/\mathbb{R}\bigr)$ is a subspace of
$l^p(G)^d$ properly containing the finite codimensional closed
subspace $\theta' (\overline{\mathbb{R} \oplus
\mathbb{R}G}_{D_p}/\mathbb{R})$.  Hence $\theta'\bigl((BD_p(G) +
\overline{\mathbb{R}G}_{D_p})/\mathbb{R} \bigr)$ is closed in
$l^p(G)^d$.  Thus with $B = \theta'\bigl((BD_p(G) +
\overline{\mathbb{R}G}_{D_p})/\mathbb{R})$ and
$A = \theta'(\overline{\mathbb{R} \oplus
\mathbb{R}G}_{D_p}/\mathbb{R})$, we obtain a contradiction from
Corollary \ref{Cmain2}.  So it is impossible for $1 <
|\partial_p(G)| < \infty$.  Therefore $|\partial_p(G)|$ is either
0, 1 or $\infty$.
\end{proof}

We will now use Corollary \ref{Cmain2} to obtain some results about
translation invariant linear functionals (which we define below) on
$D_p(G)/\mathbb{R}$.  Let $X$ be a normed space of functions on $G$.
For $f \in X$ and $h\in G$, the left translation of $f$ by $h$,
denoted by $f_h$, is the function $f_h(g) := f(hg)$.  Assume that if
$f\in X$, then $f_h \in X$ for all $h\in G$; that is, $X$ is left
translation invariant.  We shall say that $T$ is a translation
invariant left functional (TILF) on $X$ if $T(f_h) = T(f)$ for
$f \in X$ and all $h \in G$.  For the rest of this section translation
invariant will mean left translation invariant.  A common question to
ask is that if $T$ is a TILF on $X$, then is $T$ continuous?  For
background about the problem of automatic continuity, see
\cite{Meisters83, Saeki84, Willis88, Woodward74}.  Define
\[
\Diff(X) := \text{linear span}\{f_h - f \mid f\in X,\ h\in G\}.
\]
It is clear that $\Diff(X)$ is contained in the kernel of any TILF on
$X$.  Observe that $f \in D_p(G)$ if and only if $f_h - f \in l^p(G)$
for all $h \in G$.  By definition we have the following inclusions:
\[
\Diff(l^p(G)) \subseteq \Diff(D_p(G)/\mathbb{R}) \subseteq l^p(G)
\subseteq D_p(G)/\mathbb{R}.
\]
The set $D_p(G)/\mathbb{R}$ is a Banach space under the norm induced
from the norm on $D_p(G)$.  The norm of $f \in D_p(G)/\mathbb{R}$
will be indicated by $\|f \|_{D(p)}$ and the closure of a set $Y
\subseteq D_p(G)/\mathbb{R}$ will be denoted by
$\overline{Y}_{D(p)}$.  We can now state
\begin{Lem} \label{LDiff}
$\overline{\Diff(D_p(G)/\mathbb{R})}_{D(p)} =
\overline{l^p(G)}_{D(p)}$.
\end{Lem}
\begin{proof}
Let $f \in l^p(G)$.  By \cite[Lemma 1]{Woodward74} there is a
sequence $(f_n)$ in $\Diff(l^p(G))$ that converges to $f$ in the
$l^p$-norm.  It now follows from Minkowski's inequality that for $s
\in S$,
\[
\|(f-f_n)_s - (f-f_n)\|_p^p = \sum_{g\in G} |f(sg) - f_n(sg) - (f(g)
- f_n(g)) |^p \to 0
\]
as $n \to \infty$.  Hence $f \in \overline{\Diff(l^p(G))}_{D(p)}$
which implies $l^p(G) \subseteq \overline{\Diff(l^p(G))}_{D(p)}$,
and the result follows.
\end{proof}

Combining the Hahn-Banach theorem with Lemma \ref{LDiff} and the
fact $\overline{l^p(G)}_{D(p)} = D_p(G)/\mathbb{R}$ if and only if
$D_p(G)/\overline{l^p(G) \oplus \mathbb{R}}_{D_p} = 0$, we obtain the
following
\begin{Cor}
Let $1 < p \in \mathbb{R}$.  Then $D_p(G)/\overline{l^p(G) \oplus
\mathbb{R}}_{D_p} \ne 0$ if and only if there exists a nonzero
continuous TILF on $D_p(G)/\mathbb{R}$.
\end{Cor}

It was shown in \cite{Willis86} that if $G$ is nonamenable, then the
only TILF on $l^p(G)$ is the zero functional.  (Consequently every
TILF is automatically continuous!)  We now show that this is not true
for $D_p(G)/\mathbb{R}$.
\begin{Prop} \label{Ptilf}
Let $1 < p \in \mathbb{R}$ and suppose that $G$ is a
finitely generated nonamenable
group.  If $D_p(G)/\overline{l^p(G) \oplus \mathbb{R}}_{(D_p)} \ne
0$, then there is a discontinuous TILF on $D_p(G)/\mathbb{R}$.
\end{Prop}
\begin{proof}
It is known that $l^p(G)$ is closed in $D_p(G)/\mathbb{R}$ if and
only if $G$ is nonamenable \cite[Corollary 1]{Guichardet77}.  Also
$D_p(G)/(l^p(G) \oplus \mathbb{R})$ is infinite dimensional by
Corollary \ref{Cmain2}.  Let $B$ be a Hamel basis for $l^p(G)$ and
extend it to a Hamel basis $H$ for $D_p(G)/\mathbb{R}$.
Now $H \setminus B$ corresponds to a Hamel basis of $D_p(G)/(l^p(G)
\oplus \mathbb{R})$.  Select a countable subset $C := \{f_n \mid n
\in \mathbb{N}\}$ from $H\setminus B$.
Define a linear functional $T$ on $D_p(G)/\mathbb{R}$ by $T(f) = 0$
for $f \in H\setminus C$ and $T(f_n) = n\|f_n\|_{D(p)}$ for $n\in
\mathbb{N}$.  By \cite[Lemma
1]{Woodward74} $\Diff(l^p(G)) = l^p(G)$, which implies
$\Diff(D_p(G)/\mathbb{R}) = l^p(G)$.  Thus $T$ is a TILF on
$D_p(G)/\mathbb{R}$.  However
$\left(\frac{f_n}{n\|f_n\|_{D(p)}}\right) \to 0$ in
$D_p(G)/\mathbb{R}$ and $T\left(\frac{f_n}{n\|f_n\|_{D(p)}}\right) =
1$ for all $n$.  Thus $T$ is discontinuous on $D_p(G)/\mathbb{R}$.
\end{proof}

The free group on two generators provides an example of a group that
satisfies Proposition \ref{Ptilf}, see \cite[Corollary 4.3]{Puls06}
for details.  If $G$ is an infinite amenable group, then by using an
argument similar to the proof of Proposition \ref{Ptilf}, we see that
there always exists a discontinuous TILF on $D_p(G)/\mathbb{R}$ (the
key point is that $D_p(G)/(l^p(G) \oplus \mathbb{R})$ will still be
infinite dimensional).

\bibliographystyle{plain}

\end{document}